\documentclass{amsart}

\usepackage{pgf}
\usepackage{xxcolor}

\usepackage[pdftex, colorlinks=true, linkcolor=blue, urlcolor=blue, citecolor   = green]{hyperref}

\usepackage{amsmath,amssymb,amsbsy, amsthm}

\usepackage{xcolor}
\usepackage{fancyhdr, blindtext}

\usepackage{t1enc}

\usepackage{cite}

\usepackage[ msc-links]{amsrefs}

\usepackage{multirow}

 \newtheorem{thm}{Theorem}[section]
 
 \newtheorem{lem}[thm]{Lemma}
 
 \theoremstyle{definition}
 \newtheorem{defn}[thm]{Definition}
 \theoremstyle{remark}
 \newtheorem{rem}[thm]{Remark}
 \newtheorem*{ex}{Example}
 \numberwithin{equation}{section}

\def\Z{\mathbb Z}

\def\C{\mathbb C}
\def\P{\mathbb P^1 (k)}

\def\Z{{\mathbb Z}}
\def\C{{\mathbb C}}
\def\H{\mathcal H}
\def\M{\mathcal M}
\def\X{\mathcal X}
\def\S{\mathcal S}

\def\p{\mathfrak p}

\def\v{\mathfrak v}

\def\p{\mathfrak p}
\def\s{\mathfrak s}
\def\v{\mathfrak v}

\def\Jac{\mbox{Jac }}

\def\D{\Delta}

\def\a{\alpha}
\def\e{\varepsilon}
\def\b{\beta}

\def\t{\tau}
\def\<{\langle}
\def\>{\rangle}
\def\zz{\zeta}
\def\z{\omega} 

\def\iso{{\, \cong\, }}

\def\bG{\bar G}
\def\G{G}

\def\bG{{\bar G}}

\def\g{{g}}

\def\sem{{\rtimes}}

\def\Z{\mathbb Z}

\def\bC{\mathbb C}
\def\C{\mathbb C}

\def\C{\mathcal C}
\def\H{\mathcal H}
\def\M{\mathcal M}

\def\a{\alpha}
\def\b{\beta}
\def\p{\mathfrak p}

\def\Z{\mathbb Z}

\def\C{\mathbb C}

\def\H{\mathcal H}
\def\M{\mathcal M}

\def\s{\sigma}
\def\a{\alpha}
\def\p{\mathfrak p}

\def\v{\mathfrak v}
\def\<{\langle}
\def\>{\rangle}
\def\X{\mathcal X}

\def\D{\Delta}

\def\S{\mathcal S}

\def\Y{\mathcal Y}
\def\A{\mathcal A}

\def\P{\mathbb P^1 (k)}

\def\Z{{\mathbb Z}}
\def\C{{\mathbb C}}
\def\M{{\mathcal M}}
\def\H{\mathcal H}

\def\J{\mbox{Jac }}
\def\O{\Omega}

\def\Aut{\mbox{Aut }}
\def\bAut{\overline {\mathrm{Aut}}}

\def\Jac{\mbox{Jac }}

\def\embd{\hookrightarrow}

\def\D{\Delta}
\def\e{\varepsilon}

\def\g{\gamma}
\def\bG{\bar G}
\def\t{\tau}
\def\d{{\delta }}
\def\a{{\alpha }}
\def\b{{\beta }}
\def\c{\textbf{c}}

\def\bC{\mathfrak \si}

\def\<{\langle}
\def\>{\rangle}

\def\sem{{\rtimes}}

\def\s{\mathfrak s}

\def\v{\mathfrak v}

\def\p{\mathfrak p}

\def\normal{\triangleleft}

\def\zz{\zeta}

\def\z{\omega}

\def\cO{\mathcal O}

\def\Z{\mathbb Z}

\def\C{\mathbb C}

\def\H{\mathcal H}
\def\M{\mathcal M}

\def\s{\sigma}
\def\G{\Gamma}
\def\a{\alpha}
\def\b{\beta}
\def\p{\mathfrak p}
\def\P{\mathcal P}
\def\e{\varepsilon}
\def\iso{\equiv}
\def\sem{{\rtimes}}

\def\bG{\overline G}
\def\g{\gamma}

\def\iso{{\, \cong\, }}

\def\<{\langle}
\def\>{\rangle}

\def\rank{\mbox{rank }}

\title{Classifying families of superelliptic curves}

\author{Rezart Muço}
\address[Rezart Muço]{Research Institute of Science and Technology \\   Vlora, Albania}
\email{rmuco@risat.org}

\author{Nejme Pjero}
\address[Nejme Pjero]{Research Institute of Science and Technology \\ Vlora, Albania}
\email{npjero@risat.org}

\author{Ervin Ruci}
\address[Ervin Ruci]{Geolitica Inc. \\ Ottawa, Canada}
\email{eruci@risat.org}
\author{Eustrat Zhupa}
\address[Eustrat Zhupa]{UIST  "St. Paul the Apostle" \\ Ohrid, Macedonia}
\email{eustrat.zhupa@uist.edu.mk }



\begin{document}
 
\maketitle

\begin{abstract}
This paper is the first version of a project of classifying all superelliptic curves of genus $g \leq 48$ according to their automorphism group.  We determine the parametric equations in each family, the corresponding signature of the group, the dimension of the family, and the inclussion among such families.  At a later stage it will be determined the decomposition of the Jacobians and each locus in the moduli spaces of curves. 
\end{abstract}

\section{Introduction}
%*********************************************************************

In this paper we study some very classical problems  related to algebraic curves and the possibility of organizing such results in a database of curves.  Such information would be quite useful to researchers  working on coding theory, cryptography, mathematical physics, quantum computing, etc. 

In section 2 we give a brief review of the background on algebraic curves.  Throughout this paper, by "curve" we mean a smooth, irreducible algebraic curve, defined over an algebraically closed field $k$ of characteristic zero. For the reader who is interested in details and full treatment of the subject we suggest the classical books \cite{fulton, cornalba1, cornalba2, dolga_book}. 

Identifying isomorphic classes of algebraic curves is a fundamental problem of algebraic geometry which goes bach to the XIX century mathematics as a branch of invariant theory.  Especially of interest are invariants of binary forms as they determine the isomorphic classes of hyperelliptic curves and superelliptic curves (cf. Section  5).

In section 3 we give the basic background of the binary forms and their invariants.   In section 5 we discuss in more detail superelliptic curves. We define such curves and describe their automorphism groups and  the signature $\sigma$  of the covering $\X_g \to \X_g/G$.  The Hurwitz space $\H (g, G, \sigma)$ is a quasi-projective variety and  the dimension and  irreducibility can be determined as explained in section 5.  Such information enables us to fully understand such families of superelliptic curves. Identifying the isomorphism classes of such curves can be done via the invariants of binary forms (see Section 3) or the dihedral invariants as defined in \cite{super1}.

In section 5 we propose the idea of creating a database for all the algebraic curves.  The superelliptic curves consist of the bulk of the cases.  We propose what invariants need to be included in this database based on the fact that how easy it is to compute such invariants.  The main criteria of this proposed database is the automorphism group.  Each family of curves with fixed genus $g\geq 2$, group $G$, and signature $\sigma$ defines a locus in $\M_g$.  We compute the dimension and the irreducibility of such loci and an equation of such loci whenever possible.  For each family we discuss the decomposition of the corresponding Jacobians. 
Some computational packages are given as examples for curves of small genus.  We don't describe in details the mathematics behind such computer packages but for the interested reader they can be found in \cite{b-th, hyp_mod_3} and others.

%***********************************************************************
%\newpage
\section{Preliminaries on  curves}

By a \textit{curve} we mean a complete reduced algebraic curve over $\C$ which might be singular or reducible.  A \textit{smooth curve} is implicitly assumed to be irreducible. The basic invariant of a smooth curve $C$ is its genus which is half of the first Betti number of the underlying topological space.  We will denote the genus of $C$ by $g(C) = \frac 1 2 \, \rank ( H^1 (C, \Z) )$.

Let $f: \X \to \Y$ be a non-contant holomorphic map between smooth curves $\X$ and $\Y$ of genera $g$ and $g^\prime$.  For any $q\in \X$ and $p=f(q)$ in $\Y$  chose local coordinates $z$ and $w$ centered at $q$ and $p$ such that $f$ has the standart form $w= z^{\nu (q)}$. Then, for any $p$ on $C^\prime$ define \[ f^\star (p) = \sum_{q \in f^{-1} (p) } \nu ( q ) q.\]
If $D$ is any divisor on $\X$ then define $f^\star (D)$ to be the divisor on $\X$ by extending the above $f^\star (p)$. 
The degree of $n$ the divisor $f^\star (p)$ is independent of $p$ and is called the \textit{degree} of the map $f$.  The \textit{ ramification divisor } $R$ on $C$ of the map $f$ is defined by
\[ R = \sum_{q\in C} \left( \nu (p) -1 \right) q \]
The integer $\nu (q) -1$ is called the \textit{ramification index} of $f$ at $q$.   For any meromorphic differential $\phi$ on $C^\prime$ we have \[ \left( f^\star (\phi) \right) = f^\star \left( ( \phi) \right) + R\] Counting degrees we get the Riemann-Hurwitz formula 
\[ 2g-2 = n (2 g^\prime  -2) + \deg R\]
Let $\X_g$ be a genus $g\geq 2$ curve and $G$ its automorphism group (i.e, the group of automorphisms of the function field $\C (\X_g)$).  That   G is finite will be shown in section \ref{w-pts} using Weierstrass points. 

Assume $|G|=n$.  Let $L$ be the fixed subfield of $\C (\X_g)$. The field extension $\C (\X_g)/L$ correspond to a finite morphism of curves $f : \X_g \to \X_g/G$ of degree $n$.   
Denote the genus of the quotient curve $\X_g/G$ by $g^\prime$ and $R$ the ramification divisor. Assume that the covering has $s$ branch points.  Each branch point $q$ has $n/e_P$ points in its fiber $f^{-1} (q)$, where $e_P$ is the ramification index of such points $P \in  f^{-1} (q)$.

Then, $R = \sum_{i=1}^s \frac n {e_P} \left( e_P - 1 \right)$. 
By the Riemann-Hurwitz formula we have 
\[  \frac  2 n \, (g-1) = 2 g^\prime - 2 + \frac 1 n \,  \deg R  = 2g^\prime - 2 + \sum_{i=1}^s \left( 1 - \frac 1 {e_P} \right) \] 
Since  $g \geq 2$ then the left hand side is $> 0$.  Then 
\[ 2 g^\prime - 2 +  \sum_{i=1}^s \left( 1 - \frac 1 {e_P} \right) \geq 0. \] The fact that $g^\prime$, $s$, and $e_P$ are non-negative integers implies that the minimum value of this expression is 1/42.  This implies that $n \leq 84 (g-1)$. 

% Jacobians

Next we define another important invariant of the algebraic curves. Let $w_1, \dots , w_g$ be a basis of $H^0 (C, K)$ and $\gamma_1, \dots , \gamma_{2g}$ a basis for $H_1 (C, \z)$. The \textit{period matrix} $\Omega$ is the $g \times 2g$ matrix $ \Omega = \left[   \int_{\gamma_i} w_j \right]$.  The collumn vectors of the period matrix generate a lattice $\Lambda$ in $\C^g$, so that the quotient $\C^g/\Lambda$ is a complex  torus. This complex torus is called the \textit{Jacobian variety} of $C$ and denoted by the symbol $J (C)$. For more details see \cite{cornalba1} or \cite{fulton}.

If $\Lambda \subset \Lambda^\prime$ are lattices of rank $2n$ in $\C^n$ and $A= \C^n/\Lambda$, $B=\C^n/\Lambda^\prime$, then the induced map $A \to B$ is an \textit{isogeny} and $A$ and $B$ are said to be \textit{isogenous}.

% moduli space

The \textit{moduli space} $\M_g$ of curves of genus $g$ is the set of isomorphism classes of smooth, genus $g$ curves. $\M_g$ has a natural structure of a quasi-projective normal variety of dimension $3g-3$.  The compactification $\bar{\M_g}$ of $\M_g$ consists of isomorphism classes of stable curves.  A \textit{stable curve} is a curve whose only singularities are nodes and whose smooth rational components contain at least three singular points of the curve; see \cite[pg. 29]{cornalba1} and \cite[Chapter XII]{cornalba2}. $\bar{\M_g}$ is a projective variety 

Both $\M_g$ and $\bar{\M_g}$ are singular.  All the singularities arise from curves  with non-trivial automorphism group. It is precisely such curves that we intend to classify in this paper.

%********************************************************************************************************
\subsection{Automorphism groups}\label{aut-groups}

%&&&&&&&&&&&&&&&&&&&&&&&&&&&&&&&&&&&&&&&&&&&&&&&&&&&&&&&

\subsection{Superelliptic curves}\label{super}
%%%%%%%%%%%%%%%%%%%%%%%%%%%%%%%%%%%%%%%%%%%%%%%%%%%%%%%%
\def\P{\mathbb P^1}

A curve $\X$ is called superelliptic if there exist an element $\tau \in \Aut( \X)$ such that $\tau$ is central and $g \left(\X / \< \tau \> \right) =0$.  Denote by $K$ the function field of $\X_g$ and assume that the affine equation of $\X_g$ is given some polynomial in terms of $x$ and $y$. 

Let $H=\< \tau \>$ be a cyclic subgroup of $G$ such that $| H | = n$ and $H \normal G$, where $n \geq 2$. Moreover, we assume that the quotient curve $\X_g / H$ has genus zero.   The \textbf{reduced automorphism group of $\X_g$ with respect to $H$} is called the group  $\G \, := \, G/H$, see \cite{super1}, \cite{Sa1}.  

Assume  $k(x)$ is the  genus zero subfield of $K$ fixed by $H$.   Hence, $[ K : k(x)]=n$. Then, the group  $\G$ is a subgroup of the group of automorphisms of a genus zero field.   Hence, $\G <  PGL_2(k)$ and $\G$ is finite. It is a classical result that every finite subgroup of $PGL_2 (k)$  is  isomorphic to one of the following: $C_m $, $ D_m $, $A_4$, $S_4$, $A_5$.

The group $\G$ acts on $k(x)$ via the natural way. The fixed field of this action is a genus 0 field, say $k(z)$. Thus, $z$
is a degree $|\G| := m$ rational function in $x$, say $z=\phi(x)$.   $G$ is a degree $n$ extension of $\G$ and $\G$ is a finite subgroup of $PGL_2(k)$.  Hence, if we know all the possible groups that occur as $\G$ then we can determine  $G$ and the equation for $K$. The list of all groups of superelliptic curves and their equations are determined in \cite{Sa1} and \cite{Sa2}.  

Let $C$ be a superelliptic curve given by the equation \[ y^n = f(x), \] where $\deg f = d$ and $\D (f, x) \neq 0$.  
Assume that $d> n$.  Then $C$ has genus \[ g = \frac 1 2 \left( n(d-1) -d - \gcd (n, d) \frac {} {}   \right) + 1\]
Moreover, if $n$ and $d$ are relatively prime then \[ g = \frac { (n-1) (d-1) } 2,\]
see \cite{super1} for details.

%**********************************
\section{Invariants of binary forms}

In this section we define the action of $ GL_2(k)$ on binary forms and discuss the basic notions of their
invariants. Let $k[X,Z]$  be the  polynomial ring in  two variables and  let $V_d$ denote  the
$(d+1)$-dimensional  subspace  of  $k[X,Z]$  consisting of homogeneous polynomials.
\begin{equation}  \label{eq1}
f(X,Z) = a_0X^d + a_1X^{d-1}Z + \dots + a_dZ^d
\end{equation}
of  degree $d$. Elements  in $V_d$  are called  {\it binary  forms} of degree $d$.  We let $GL_2(k)$ act as a
group of automorphisms on $ k[X, Z] $   as follows:
\begin{equation}
 M =
\begin{pmatrix} a &b \\  c & d
\end{pmatrix}
\in GL_2(k), \textit{   then       }
\quad  M  \begin{pmatrix} X\\ Z \end{pmatrix} =
\begin{pmatrix} aX+bZ\\ cX+dZ \end{pmatrix}
\end{equation}
This action of $GL_2(k)$  leaves $V_d$ invariant and acts irreducibly on $V_d$.
\begin{rem}
It is well  known that $SL_2(k)$ leaves a bilinear  form (unique up to scalar multiples) on $V_d$ invariant. This
form is symmetric if $d$ is even and skew symmetric if $d$ is odd.
\end{rem}
Let $A_0$, $A_1$,   \dots  , $A_d$ be coordinate  functions on $V_d$. Then the coordinate  ring of $V_d$ can be
identified with $ k[A_0  ,  \dots  , A_d] $. For $I \in k[A_0,  \dots  , A_d]$ and $M \in GL_2(k)$, define $I^M \in k[A_0,
 \dots  ,A_d]$ as follows
\begin{equation} \label{eq_I}
{I^M}(f):= I(M(f))
\end{equation}
for all $f \in V_d$. Then  $I^{MN} = (I^{M})^{N}$ and Eq.~(\ref{eq_I}) defines an action of $GL_2(k)$ on $k[A_0,
 \dots  ,A_d]$.
A homogeneous polynomial $I\in k[A_0, \dots , A_d, X, Z]$ is called a {\it covariant}  of index $s$ if
$$I^M(f)=\delta^s I(f),$$
where $\delta =\det(M)$.  The homogeneous degree in $a_1, \dots , a_n$ is called the {\it degree} of $I$,  and the
homogeneous degree in $X, Z$ is called the {\it  order} of $I$.  A covariant of order zero is called {\it
invariant}.  An invariant is a $SL_2(k)$-invariant on $V_d$.

We will use the symbolic method of classical theory to construct covariants of binary forms.    Let
$$f(X,Z):=\sum_{i=0}^n
\begin{pmatrix} n \\ i
\end{pmatrix}
a_i X^{n-i} \, Z^i, \quad  and \quad g(X,Z) :=\sum_{i=0}^m
  \begin{pmatrix} m \\ i
\end{pmatrix}
b_i X^{n-i} \, Z^i
$$
be binary forms of  degree $n$ and $m$ respectively with coefficients in $k$. We define the {\bf r-transvection}
$$(f,g)^r:= \frac {(m-r)! \, (n-r)!} {n! \, m!} \, \, \sum_{k=0}^r
(-1)^k
\begin{pmatrix} r \\ k
\end{pmatrix} \cdot
\frac {\partial^r f} {\partial X^{r-k} \, \,  \partial Z^k} \cdot \frac {\partial^r g} {\partial X^k  \, \,
\partial Z^{r-k} }
 $$
 It is a homogeneous polynomial in $k[X, Z]$ and therefore a covariant
of order $m+n-2r$ and degree 2. In general, the $r$-transvection of two covariants of order $m, n$ (resp., degree
$p, q$) is a covariant of order $m+n-2r$  (resp., degree $p+q$).

For the rest of this paper $F(X,Z)$ denotes a binary form of order $d:=2g+2$ as below
\begin{equation}
F(X,Z) =   \sum_{i=0}^d  a_i X^i Z^{d-i} = \sum_{i=0}^d
\begin{pmatrix} n \\ i
\end{pmatrix}    b_i X^i Z^{n-i}
\end{equation}
where $b_i=\frac {(n-i)! \, \, i!} {n!} \cdot a_i$,  for $i=0, \dots , d$.  We denote invariants (resp.,
covariants) of binary forms by $I_s$ (resp., $J_s$) where the subscript $s$ denotes the degree (resp., the order).
We define the following covariants and invariants:
\begin{equation}
\begin{split}\label{covar}
\aligned
 I_2 & :=(F,F)^d,   \\
 I_4 & :=(J_4, J_4)^4,  \\
 I_6 & :=((F, J_4)^4, (F, J_4)^4)^{d-4},   \\
 I_6^\ast & :=((F, J_{12})^{12}, (F, J_{12})^{12})^{d-12},  \\
  M  & :=((F, J_4)^4, (F, J_8)^8)^{d-10}, \\
\endaligned
\qquad \aligned
& J_{4j}   :=   (F,F)^{d-2j}, \, \,  j=1, \dots , g, \\
& I_4'    :=     (J_8, J_8)^8, \\
& I_6^\prime  :=((F, J_8)^8, (F, J_8)^8)^{d-8}, \\
& I_3      :=(F, J_d)^d, \\
& I_{12}   :=(M, M)^8\\
\endaligned
\end{split}
\end{equation}
{\it Absolute invariants} are called $GL_2(k)$-invariants. We  define the following absolute invariants:
$$i_1:=\frac {I_4'} {I_2^2}, \, \,  i_2:=\frac {I_3^2} {I_2^3},
\, \,  i_3:=\frac {I_6^\ast  } {I_2^3}, \, \,  j_1 := \frac {I_6^{'}} {I_3^2}, \, \,   j_2:= \frac {I_6} {I_3^2},
\, \, s_1:=\frac {I_6^2} {I_{12}}, \, \, s_2:=\frac {(I_6^{'})^2} {I_{12}}
$$
$$
\v_1:= \frac {I_6} {I_6^\ast }, \, \,  \v_2:=\frac {(I_4^{'})^3} {I_3^4}, \, \, \v_3:= \frac {I_6} {I_6^{'}}, \,
\,  \v_4:=\frac {(I_6^\ast )^2} {I_4^3}.
$$
In the case $g=10$ and $I_{12}=0$ we define
\begin{equation}
\begin{split}
 I_6^\star & := ( (F, J_{16})^{16},  (F, J_{16})^{16} )^{d-16}), \\
S   \, \, & :=( J_{12}, J_{16} )^{12}, \\
I_{12}^\ast & := ( \, (J_{16}, S)^4, \, (J_{16}, S)^4 \, )^{12}\\
\end{split}
\end{equation}
and $$\v_5\,  :=\frac {I_6^\star }{I_{12}^\ast}.$$

\noindent For a given curve $\X_g$ we denote by $I(\X_g)$ or $i(\X_g)$ the corresponding  invariants.

\subsection{Invariants of binary Sextics}
%**************************************
Let $f(X,Y)$ be the binary sextic
$$f(X,Y) = \sum_{i=0}^6 a_i X^i Y^{6-i}.$$
defined over an algebraically closed field $k$.  We define the following covariants:
\begin{equation}
\begin{split}
&H=(f,f)^2, \quad i=(f, f )^4, \quad l=(i,f)^3, \\
\end{split}
\end{equation}
Then, the following 
\begin{equation}\label{def_J}
\begin{aligned}
&     J_2= (f,f)^6,                    &   \qquad   &  J_4 = (i,i )^4, \\
&  J_6=  (l,l)^2,                      &    \qquad      &   J_{10}=(f,l^3)^6,  \\
\end{aligned}
\end{equation}
are $SL_2(k)$- invariants; see \cite{hyp_mod_3} for details.     

The  	absolute invariants $t_1,t_2$ and $t_3$ called absolute invariants are:
\begin{equation}\label{def_Igusa}
\begin{aligned}
&     t_1= \frac{J_2^5}{J_{10}},  &   \qquad   &  t_2 = \frac{J_2^3 \cdot J_4}{J_{10}}, &   \qquad   &  t_3 = \frac{J_2^2 \cdot J_6}{J_{10}} , \\
\end{aligned}
\end{equation}

\begin{lem}
Two genus two curves $C$ and $C^\prime$ are isomorphicm if and only if they have the same absolute invariants.
\end{lem}

\subsection{Invariants of binary octavics}
%**************************************

Let $f(X,Y)$ be the binary octavic
$$f(X,Y) = \sum_{i=0}^8 a_i X^i Y^{8-i}.$$
defined over an algebraically closed field $k$.  We define the following covariants:
\begin{equation}
\begin{split}
&g=(f,f)^4, \quad k=(f, f )^6, \quad h=(k,k)^2, \\
& m=(f,k)^4, \quad  n=(f,h)^4, \quad p=(g,k)^4, \quad q=(g, h)^4,\\
\end{split}
\end{equation}
where the operator $( \cdot    , \cdot )^n$ denotes the $n$-th transvection of two binary forms; see \cite{hyp_mod_3} among many other references.   Then, the following 
\begin{equation}\label{def_J}
\begin{aligned}
&     J_2= 2^2 \cdot 5 \cdot 7 \cdot (f,f)^8,                    &   \qquad   &  J_3= \frac 1 3 \cdot  2^4 \cdot 5^2 \cdot 7^3 \cdot (f,g )^8, \\
&  J_4= 2^9 \cdot 3 \cdot 7^4 \cdot (k,k)^4,                      &    \qquad      &   J_5= 2^9 \cdot 5 \cdot 7^5 \cdot (m,k)^4,  \\
&   J_6 = 2^{14} \cdot 3^2 \cdot 7^6 \cdot (k,h )^4,              &    \qquad      &      J_7= 2^{14} \cdot 3 \cdot 5 \cdot 7^7 \cdot (m,h )^4,  \\
& J_8= 2^{17} \cdot 3 \cdot 5^2 \cdot 7^9 \cdot   (p,h)^4,        &    \qquad       &   J_9= 2^{19} \cdot 3^2 \cdot 5 \cdot 7^9 \cdot   (n,h)^4, \\
&  J_{10}=   2^{22} \cdot 3^2 \cdot 5^2 \cdot 7^{11} (q,h)^4      &  \qquad    &       \\
\end{aligned}
\end{equation}
are $SL_2(k)$- invariants; see \cite{hyp_mod_3} for details.     

Next, we define  $GL(2, k)$-invariants   as  follows
\[
t_1:= \frac {J_3^2} {J_2^3}, \quad 
t_2:= \frac {J_4} {J_2^2}, \quad 
t_3:= \frac {J_5} {J_2\cdot J_3}, \quad 
t_4:= \frac {J_6} {J_2\cdot  J_4}, \quad 
t_5:= \frac {J_7} {J_2 \cdot J_5}, \quad 
t_6:= \frac {J_8} {J_2^4},
\]
There is an algebraic relation 
\begin{equation}\label{main_eq}
T(t_1, \dots , t_6)=0\
\end{equation}
that such invariants satisfy, computed in \cite{hyp_mod_3}.   
The field of invariants $\S_8$ of binary octavics is $\S_8= k(t_1, \dots , t_6),$  where $t_1, \dots , t_6$ satisfies the equation $T(t_1, \dots , t_6)=0$. Hence, we have an explicit description of the hyperelliptic moduli $\H_3$; see \cite{hyp_mod_3} for details. 

Throughout this paper we will use the following important result

\begin{lem}[Shaska \cite{hyp_mod_3}]
Two genus 3 hyperelliptic curves $C$ and $C^\prime$, defined over an algebraically closed field $k$ of characteristic zero, with $J_2, J_3, J_4, J_5$ nonzero are isomorphic over $k$ if and only if \[ t_i (C) = t_i(C^\prime), \quad \textit{ for } \quad i=1, \dots 6.\]
\end{lem}

In the cases of curves when $t_1, \dots , t_6$ are not defined we will define new invariants as suggested in \cite{hyp_mod_3}.   From \cite[Lemma 4]{hyp_mod_3} we know that $J_2, \dots , J_7$ can't all be 0, otherwise the binary form would have a multiple root.

\subsection{Invariants of binary decimics}
%**************************************

Let $f(X,Y)$ be the binary decimic defined over an algebraically closed field $k$.  We define the following covariants:
\begin{equation}
\begin{split}
&k=(f,f)^8, \quad q=(f, f )^6, \quad m=(m,k)^4, \\
& r=(f,q)^8, \quad  k_q=(q,q)^6, \quad k_m=(m,m)^4, \quad m_q=(q,k_q)^4\\
\end{split}
\end{equation}
\noindent and the following invariants:
\begin{equation}
\begin{split}
\begin{aligned}
&     J_2= (f,f)^{10},                    &   \qquad   &  A_6 = (m,m)^6, \\
&  J_4=  (k,k)^4,                      &    \qquad      &   C_6=(r,r)^2,  \\
&  J_8=  (k,k_m)^4,                      &    \qquad      &   J_{14}=((k_q,k_q)^2,m_q)^4,  \\
&  J_9=  ((k,m)^1,k \cdot k)^8,                      &    \qquad      &   A_{14}=((k,k)^2 \cdot (k,k)^2,(m,m)^2)^8,  \\
&  J_{10}=  ((m,m)^2,k \cdot k)^8 \\
\end{aligned}
\end{split}
\end{equation}

\begin{thm}The eight invariants $J_2$, $J_4$, $A_6$, $C_6$, $J_8$, $J_9$, $J_{10}$, $J_{14}+A_{14}$ form a homogeneous system of parameters for the ring $\cO(V_n)^{SL_2}$ of invariants of the binary decimics.
\end{thm}

For the proof the reader can check \cite{binary-dec}. We still do not know a set of $GL_2 (k)$-invariants for the binary decimics. 

%*******************************************************
%\newpage
\section{A list of all superelliptic curves with extra automorphisms for $g\leq 10$} 

To create a database of algebraic curves which contains enough information we have to start with superelliptic curves as the simplest cases of all the families.  It turns out that the superelliptic curves are the overwhelming majority of automorphism groups of curves for any fixed genus. 

In the tables below are displayed all superelliptic curves for genus $5 \leq g \leq 10$. The cases for $g <5$ have appeared before in the literature.  The first column of the table represents the case from Table~1 of \cite{Sa1}, the second column is the reduced automorphism group.  In the third column we put information about the full automorphism group. Such groups are well known and we only display the 'obvious' cases, for full details one can check \cite{Sa1} and \cite{Sa2}.  

In the fourth column is the level $n$ of the superelliptic curve; see \cite{super1}.  Hence, the equation of the curve is given by 
\[ y^n = f(x), \]
where $f(x)$ is the polynomial displayed in the last column.  Columns 5 and 6 respectively represent the order of an automorphism in the reduced automorphism group and the signature of the covering $\X \to \X/ G$. The sixth column represents the dimension of the corresponding locus in the moduli space $\M_g$. 
Throughout these tables $f_1 (x)$ is as follows
\[ f_1(x) = x^{12}-a_1 x^{10}-33x^8+2a_1-33x^4-a_1x^2+1 \]
The signatures of the coverings are not fully given.  Indeed, for a full signature $\left( \sigma_1, \dots , \sigma_r\right)$ we know that $\sigma_1 \cdots \sigma_r=1$.  Hence, $\sigma_r = \sigma_{r-1}^{-1}\cdots \sigma_1^{-1}$.  We do not write $\sigma_r$. 

%\subsection{Superelliptic curves of genus $g=5$}

\begin{table}[hbt]
\caption{Superelliptic curves for genus $5 \leq g \leq 10$}
\begin{tabular}{|l|l|l|l|l|l|l|l|}
\hline
Nr. & $\bar G$             & G&$n$  &$m$ & sig. & $\delta$ & Equation $y^n=f(x)$ \\
\hline \hline
\multicolumn{8}{c}{Genus 5} \\
\hline \hline
 1& \multirow{4}{*}{$C_m$} &   $C_2^2$        &  2  & 2  & $2^7$      & 5 & $x^{12} + \sum_{i=1}^5 a_i x^{2i} + 1$  \\ 
 1&                        &   $C_3 \times C_2$        &  2  & 3  & $2^3, 3^2$ & 3 & $x^{12} + \sum_{i=1}^3 a_i x^{3i} + 1$ \\   
 1&                        &   $C_2\times C_4$        &  2  & 4  & $2^2, 4^2$ & 2  & $x^{12} + a_2 x^8 + a_1 x^4 +1$ \\  
 2&                        & $C_{22}$  &  2  & 11 &   11, 22   & 0   & $x^{11}+1$ \\ 
 2&                        & $C_{22}$  &  11 & 2  &   2, 22   &  0    & $x^2+1$ \\ 
 3&                        & $C_2$     &  2  & 1  &  $2^{11}$   & 9 & $x^{10}+\sum_{i=1}^9 a_ix^i +1$ \\ 
 3&                        & $C_4$     &  2  & 2  & 2         &   4 & $x^{10}+\sum_{i=1}^4 a_i x^{2i} +1$ \\ 
\hline
 4& \multirow{4}{*}{$D_{2m}$} &  &  2 & 2  & $2^6$        &  3 & $ \prod_{i=1}^3 (x^4+a_i x^2+1) $ \\ 
 4&                           &  &  2 & 3  & $2^4$, 3     &  2 & $(x^6+a_1x^3+1)(x^6+a_2x^3+1)$ \\   
 4&                           &  &  2 & 6  & $2^3$, 6     &  1 & $x^{12}+a_1x^6+1$ \\  
 5&                           &  &  2 & 4  & $2^2$, $4^2$ &  1 & $(x^4-1)(x^8+ a_1 x^4 +1)$  \\ 
 5&                           &  &  2 & 12 & 2, 4, 12     &  0 & $x^{12}-1$  \\ 
 6&                           &  &  2 & 5  & $2^3$, 10    &  1 & $x(x^{10}+a_1x^5+1)$ \\ 
 7&                           &  &  2 & 2  & $2^3$, $4^2$ &  2 & $(x^4-1)(x^4+a_1x^2+1)(x^4+a_2x^2+1)$ \\  
 7&                           &  &  2 & 3  & 2, 3, $4^2$  &  1 & $(x^6-1)(x^6+a_1x^3+1)$ \\  
 8&                           &  &  2 & 2  & $2^3$, $4^2$ &  2 & $x(x^2-1)(x^4+a_1x^2+1)(x^4+a_2x^2+1)$ \\  
 8&                           &  &  2 & 10 & 2, 4, 20     &  0 & $x (x^{10} -1)$ \\   
\hline   
10 &   $A_4$                     &  & 2 &  & $2^2$, $3^2$ & 1 & $f_1(x)$ \\  
\hline   
20 &   $S_4$                     &  & 2  & 0 & 3, $4^2$ & 0 & $x^{12} - 33x^8-33x^4+1$ \\  
\hline   
25 &   $A_5$                     &  & 2 &  & 2,3,10 & 0 & $x (x^{10} + 11x^5-1)$ \\  
\hline \hline
\multicolumn{8}{c}{Genus 6} \\
\hline \hline
 1& \multirow{4}{*}{$C_m$} &   $C_2^2$        &  2  & 2  & $2^8$      & 6 & $x^{14} + \sum_{i=1}^6 a_i x^{2i} + 1$  \\
 2&                        &   $C_{26}$        &  2  & 13  & $13, 26$ & 0 & $x^{13} + 1$ \\
 2&                        &   $C_{21}$        &  3  & 7  & $7, 21$ & 0  & $x^7 +1$ \\
 2&                        & $C_{20}$  &  4  & 5 &   $5, 20$   & 0   & $x^{5}+1$ \\
 2&                        & $C_{10}$  &  5 & 2  &   $2, 5, 10$  &  1    & $x^4+a_1x^2+1$ \\
 2&                        & $C_{20}$     &  5  & 4  &  $4, 20$   & 0 & $x^4+1$ \\
 2&                        & $C_{21}$     &  7  & 3  & $3, 21$    & 0 & $x^3 +1$ \\
 2&                        & $C_{26}$     &  13 & 2  & $2, 26$    & 0 & $x^2 +1$ \\
 3&                        & $C_{2}$     &  2  & 1  & $2^{13}$    & 11 & $x^{12} +\sum_{i=1}^{11} a_i x^{i}+1$ \\
 3&                        & $C_{4}$     &  2  & 2  & $2^5, 4^2$    & 5 & $x^{12} +\sum_{i=1}^{5} a_i x^{2i}+1$ \\
 3&                        & $C_{6}$     &  2  & 3  & $2^3, 6^2$    & 3 & $x^{12} +\sum_{i=1}^{3} a_i x^{3i}+1$ \\
 3&                        & $C_{8}$     &  2  & 4  & $2^2, 8^2$    & 2 & $x^{12}+ \sum_{i=1}^{2} a_i x^{4i}+1$ \\
 3&                        & $C_{3}$     & 3  & 1  & $3^7$    & 5 & $x^{6}+ \sum_{i=1}^{5} a_i x^{i}+1$ \\
 3&                        & $C_{6}$     &  3  & 2  & $3^2, 6^2$    & 2 & $x^{6}+a_2x^4+a_1x^2+1$ \\
 3&                        & $C_{4}$     &  4  & 1  & $4^5$    & 3 & $x^{4}+ \sum_{i=1}^{3} a_i x^{i}+1$ \\
 3&                        & $C_{5}$     &  5  & 1  & $5^4$    & 2 & $x^3 +a_1x+a_2x^2+1$ \\
\hline
 4& \multirow{4}{*}{$D_{2m}$} &$ D_{14}\times C_2$ &  2 & 7  & $2^3, 7$        &  1 & $x^{14}+a_1x^7+1)$ \\
 5&                           & $G_5$ &  2 & 2  & $2^5$, 4     &  3 & $(x^2-1) \prod_{i=1}^3 (x^4+a_i x^2+1) $ \\
 5&                           & $G_5$ &  2 & 14  & $2, 4, 14$     &  0 & $x^{14}- 1$ \\
 5&                           & $ D_{10}\times C_2$ &  5 & 5  & $2, 5, 10$ &  0 & $x^5-1$  \\
 6&                           & $D_8$ &  2 & 2 & $2^5, 4 $     & 3 & $x \cdot  \prod_{i=1}^3 (x^4+a_i x^2+1)  $  \\
 6&                           &$ D_{6}\times C_2$  &  2 & 3  & $2^4$, 6    &  2 & $x \cdot \prod_{i=1}^2 (x^4+a_i x^2+1)$ \\
\hline
\end{tabular}
\end{table}

\addtocounter{table}{-1}
\begin{table}
\caption{(Cont.)}
\begin{tabular}{|l|l|l|l|l|l|l|l|}
\hline
Nr. & $\bar G$             & G&$n$  &$m$ & sig. & $\delta$ & Equation $y^n=f(x)$ \\
\hline  
 6&                           &$D_{24}$  &  2 & 6  & $2^3$, 12 &  1 & $x(x^{12}+a_1x^6+1)$ \\
 6&                           & $ D_{6}\times C_3$ &  3 & 3  &  $2^2, 3, 9$  &  1 & $x(x^6+a_1x^3+1)$ \\
 6&                           & $ D_{16}$ &  4 & 2  & $2^2$, 4, 8 &  1 & $x(x^4+a_1x^2+1)$ \\
 8&                           & $G_8$ &  2 & 4 & $2^2$, 4, 8     &  1 & $x (x^{4} -1)(x^8+a_1x^4+1)$ \\
 8&                           & $G_8$ &  2 & 12 & 2, 4, 24     &  0 & $x (x^{12} -1)$ \\
 8&                           &$ D_{4}\times C_3$  &  3 & 2 & 2, 3, $6^2$     &  1 & $x (x^{2} -1)(x^4+a_1x^2+1)$ \\
 8&                           & $ D_{12}\times C_3$ &  3 & 6 & 2, 6, 18     &  0 & $x (x^{6} -1)$ \\
 8&                           &$G_8$  &  4 & 4 & 2, 8, 16     &  0 & $x (x^{4} -1)$ \\
 8&                           & $ D_{6}\times C_5$ &  5 &3 & 2, 10, 15     &  0 & $x (x^{3} -1)$ \\
 8&                           &$ D_{4}\times C_7$  &  7 & 2 & 2, $14^2$     &  0 & $x (x^{2} -1)$ \\
 9&                           &$G_9$  &  2 & 2 & $2^2$, $4^3$     &  2 & $x (x^{4} -1) \cdot \prod_{i=1}^2 (x^4+a_i x^2+1)$ \\
 9&                           & $G_9$ &  2 & 3 & $2, 4^2, 6$     &  1 & $x (x^{6} -1)(x^6+a_1x^3+1)$ \\
\hline
% &   $A_4$                     &  &  &  &  & &  \\
%\hline
 18&   $S_4$                     & $ G_{18}$ & 4  & 0 & 2, 3, 16 & 0 & $x(x^4-1)$ \\
 19&                             & $ G_{19}$ & 2  & 0 & 2, 6, 8 & 0 & $x(x^4-1)(x^8+14x^4+1)$ \\
\hline \hline
\multicolumn{8}{c}{Genus 7} \\
\hline \hline
 1& \multirow{4}{*}{$C_m$} &   $C_2^2$        &  2  & 2  & $2^9$      & 7 & $x^{16} + \sum_{i=1}^7 a_i x^{2i} + 1$  \\
 1&                        &   $C_2 \times C_4$        &  2  & 4 & $2^3, 4^2$ & 3 & $x^{16} + \sum_{i=1}^3 a_i x^{4i} + 1$ \\
 1&                        &   ${C_3}^2$        &  3  & 3  & $3^4$ & 2  & $x^{9} + a_2 x^6 + a_1 x^3 +1$ \\
 2&                        & $C_{6}$  &  2  & 3 &  $2^4, 3, 6$    & 4   & $x^{15}+ \sum_{i=1}^4 a_1 x^{3i} +1$ \\
 2&                        & $C_{10}$  &  2 & 5  &   $2^2, 5, 10$   &  2    & $x^{15}+a_1x^5+a_2x^{10}+1$ \\
 2&                        & $C_{30}$     &  2  & 15  &  15, 30   & 0 & $x^{15} +1$ \\
 2&                        & $C_6$     &  3  & 2  & 2, $3^3$, 6         &   3 & $x^8+a_3x^6+a_2x^4+a_1x^2 +1$ \\
 2&                        & $C_{12}$     &  3  & 4  & 3, 4, 12         &   1 & $x^8+a_1x^4+1$ \\
 2&                        & $C_{24}$     &  3  & 8  & 8, 24         &   0 & $x^8+1$ \\
 2&                        & $C_{30}$     &  15  & 2  & 2, 30         &   0 & $x^2+1$ \\
 3&                        & $C_{2}$     &  2  & 1  & $2^{15}$         &   13 & $x^{14}+\sum_{i=1}^{13} a_i x^{i}+1$ \\
 3&                        & $C_{4}$     &  2  & 2  & $2^6, 4^2$       &   6 & $x^{14}+\sum_{i=1}^6 a_i x^{2i}+1$ \\
 3&                        & $C_{3}$     &  3  & 1  & $3^8$       &   6 & $x^{7}+\sum_{i=1}^6 a_i x^{i}+1$ \\
\hline
 4& \multirow{4}{*}{$D_{2m}$} & $D_4\times C_2$ &  2 & 2  & $2^7$        &  4 & $\prod_{i=1}^4(x^4+a_ix^2+1)$ \\
 4&                           &$D_8\times C_2$  &  2 & 4  & $2^4$, 4     &  2 & $(x^8+a_1x^4+1)(x^8+a_2x^4+1)$ \\
 4&                           &$D_{16}\times C_2$  &  2 & 8  & $2^3$, 8  &  1 & $x^{16}+a_1x^8+1$ \\
 5&                           &$G_5$  &  2 & 16  & 2, 4, 16 &  0 & $x^{16}-1$  \\
 5&                           & $D_6\times C_3$ &  3 & 3 & 2, $3^2$, 6    &  1 & $(x^3 - 1)(x^6+a_1x^3+1)$  \\
 5&                           & $D_{18}\times C_3$ &  3 & 9 & 2, 6, 9  &  0 & $x^9-1$  \\
 6&                           & $D_{14}\times C_2$ &  2 & 7  & $2^3$, 14    &  1 & $x(x^{14}+a_1x^7+1)$ \\
 7&                           &$G_7$  &  2 & 2  & $2^4$, $4^2$ &  3 & $(x^4-1) \, \prod_{i=1}^3 (x^4+a_i x^2+1) $ \\
 7&                           & $G_7$ &  2 & 4  & 2, $4^3$  &  1 & $(x^8-1)(x^8+a_1x^4+1)$ \\
 8&                           & $G_8$ &  2 & 2  & $2^4$, $4^2$ &  3 & $x(x^2-1)\, \prod_{i=1}^3 (x^4+a_i x^2+1) $ \\
 8&                           & $G_8$ &  2 & 14 & 2, 4, 28     &  0 & $x (x^{14} -1)$ \\
 8&                           & $D_{14}\times C_3$ &  3 & 7 & 2, 6, 21     &  0 & $x (x^{7} -1)$ \\
 8&                           & $G_8$ &  8 & 2 & 2,${16}^2$     &  0 & $x (x^{2} -1)$ \\
\hline \hline
11 &   $A_4$                     &$K$  & 2 & 0  & $2^2$, 3, 6 & 1 & $(x^4+2 \sqrt{-3} x^2+1) \, f_1(x)$ \\
\hline \hline
\end{tabular}
\end{table}

\addtocounter{table}{-1}
\begin{table}
\caption{(Cont.)}
\begin{tabular}{|l|l|l|l|l|l|l|l|}
\hline
Nr. & $\bar G$             & G&$n$  &$m$ & sig. & $\delta$ & Equation $y^n=f(x)$ \\
\hline \hline
\multicolumn{8}{c}{Genus 8} \\
\hline \hline
 1& \multirow{4}{*}{$C_m$} &   $C_2^2$        &  2  & 2  & $2^{10}$      & 8 & $x^{18} + \sum_{i=1}^8 a_i x^{2i} + 1$  \\
 1&                        &   $C_2 \times C_3$        &  2  & 3  & $2^5, 3^2$ & 5 & $x^{18} + \sum_{i=1}^5 a_i x^{3i} + 1$ \\
 1&                        &   $C_2\times C_6$        &  2  & 6  & $2^2, 6^2$ & 2  & $x^{18} + a_1 x^6 + a_2 x^{12} +1$ \\
 2&                        & $C_{34}$  &  2  & 17 &   17, 34   & 0   & $x^{17}+1$ \\
 2&                        & $C_{34}$  &  17 & 2  &   2, 34   &  0    & $x^2+1$ \\
 3&                        & $C_2$     &  2  & 1  &  $2^{17}$   & 15 & $x^{16}+\sum_{i=1}^15 a_ix^i +1$ \\
 3&                        & $C_4$     &  2  & 2  & $2^7, 4^2$        &  7 & $x^{16}+\sum_{i=1}^7 a_i x^{2i} +1$ \\
  3&                        & $C_8$     &  2  & 4  & $2^3, 8^2$        & 3 & $x^{16}+a_1x^4+a_2x^8+a_3x^{12} +1$ \\
\hline
 4& \multirow{4}{*}{$D_{2m}$} &$D_6\times C_2$  &  2 & 3  & $2^5$, 3    & 3 & $ \prod_{i+1}^3 (x^6+a_i x^3+1) $ \\
 4&                           &$D_{18}\times C_2$  &  2 & 9  & $2^3$, 9 & 1 & $x^{18}+a_1x^9+1$ \\
 5&                           & $G_5$ &  2 & 2  & $2^6$, 4     &  4 & $(x^2-1)\prod_{i=1}^4(x^4+a_ix^2+1)$ \\
 5&                           &$G_5$  &  2 & 6  & $2^2$, 4, 6 &  1 & $(x^6-1)(x^{12}+a_1x^6+1)$  \\
 5&                           &$G_5$  &  2 & 18 & 2, 4, 18     &  0 & $x^{18}-1$  \\
 6&                           &$D_8$  &  2 & 2  & $2^6$, 4    &  4 & $x\prod_{i=1}^4(x^4+a_ix^2+1)$ \\
 6&                           &$D_{16}$  &  2 & 4  & $2^4$, 8    &  2 & $x(x^8+a_1x^4+1)(x^8+a_2x^4+1)$ \\
 6&                           &$D_{32}$  &  2 & 8  & $2^3$, 16    &  1 & $x(x^{16}+a_1x^8+1)$ \\
 7&                           &$G_9$  &  2 & 3  & $2^2$, 3,  $4^2$ &  2 & $(x^6-1)(x^6+a_1x^3+1)(x^6+a_2x^3+1)$ \\
 8&                           &$G_8$  &  2 & 16  & 2, 4, 32 &  0 & $x(x^{16}-1)$ \\
 9&                           & $G_9$ &  2 & 2 & $2^3, 4^3$    &  3 & $x(x^4-1)\, \prod_{i+1}^3 (x^6+a_i x^3+1) $ \\
  9&                           & $G_9$ &  2 & 4 & 2, $ 4^2$, 8   &  1 & $x(x^8-1)(x^8+a_1x^4+1)$ \\
\hline
13 &   $A_4$                     & $K$ & 2 & 0 & 2, $3^2$, 4 & 1 & $x(x^4-1)\, f_1 (x)$ \\
\hline
22 &   $S_4$                     &$G_{22}$  & 2  & 0 & 3, 4, 8 & 0 & $x(x^4-1)(x^{12}-33x^8-33x^4+1)$ \\
\hline \hline
\multicolumn{8}{c}{Genus 9} \\
\hline \hline
 1& \multirow{4}{*}{$C_m$} &   $C_2^2$        &  2  & 2  & $2^{11}$      & 9 & $x^{20} + \sum_{i=1}^9 a_i x^{2i} + 1$  \\
 1&                        &   $C_2 \times C_4$ &  2  & 4  & $2^4, 4^2$ & 4 & $x^{20} + \sum_{i=1}^4 a_i x^{4i} + 1$ \\
 1&                        &   $C_2\times C_5$  &  2  & 5 & $2^3, 5^2$ & 3  & $x^{20} + a_1 x^5 + a_2 x^{10}+a_3x^{15} +1$ \\
 1&                        &   $C_2\times C_4$  &  4  & 2 & $2^2, 4^3$ & 3  & $x^{8} + a_1 x^2 + a_2 x^{4}+a_3x^{6} +1$ \\
 2&                        & $C_{38}$  &  2  & 19 &   19, 38   & 0   & $x^{19}+1$ \\
 2&                        & $C_{6}$  &  3 & 2  &   2, $3^4$, 6   &  4    & $x^{10} + a_1 x^2 + a_2 x^{4}+a_3x^{6}+a_4x^8 +1$ \\
 2&                        & $C_{15}$  &  3 & 5  &   3, 5, 15  &  1   & $x^{10} + a_1 x^5 +1$ \\
 2&                        & $C_{30}$  &  3 & 10  &  10, 30  &  0    & $x^{10} +1$ \\
 2&                        & $C_{28}$  &  4 & 7  &  7, 28   &  0    & $x^7 +1$ \\
 2&                        & $C_{14}$  &  7 & 2  &   2,7, 14   &  1   & $x^4+a_1x^2+1$ \\
  2&                        & $C_{28}$  &  7 & 4  &  4, 28   & 0    & $x^4 +1$ \\
  2&                        & $C_{30}$  &  10 & 3  &  3, 30    &  0      & $x^3 +1$ \\
  2&                        & $C_{38}$  &  19 & 2  &   2, 38   &  0      & $x^2  +1$ \\

 3&                        & $C_2$     &  2  & 1  &  $2^{19}$   & 17 & $x^{18}+\sum_{i=1}^{17} a_ix^i +1$ \\
 3&                        & $C_4$     &  2  & 2  & $2^8, 4^2$        &   8 & $x^{18}+\sum_{i=1}^8 a_i x^{2i} +1$ \\
 3&                        & $C_6$     &  2  & 3  & $2^5, 6^2$        &   5 & $x^{18}+\sum_{i=1}^5 a_i x^{3i} +1$ \\
 3&                        & $C_{12}$     &  2  & 6  & $2^2, {12}^2$  &   2 & $x^{18}+ a_1x^6+a_2x^{12} +1$ \\
 3&                        & $C_{3}$     &  3   & 1  & $3^{10}$  &   8 & $x^{9}+ \sum_{i=1}^8 a_i x^{i}+1$ \\
 3&                        & $C_{9}$     &  3   & 3  & $3^2, 9^2$  &   2 & $x^{9}+ a_2x^6+a_1x^3+1$ \\
\hline \hline
\end{tabular}
\end{table}

\addtocounter{table}{-1}
\begin{table}
\caption{(Cont.)}
\begin{tabular}{|l|l|l|l|l|l|l|l|}
\hline 
\hline
Nr. & $\bar G$             & G&$n$  &$m$ & sig. & $\delta$ & Equation $y^n=f(x)$ \\
\hline 
 3&                        & $C_{4}$     &  4   & 1  & $4^7$  &   5 & $x^{6}+ \sum_{i=1}^5 a_i x^{i}+1$ \\
 3&                        & $C_{8}$     &  4   & 2  & $4^2, 8^2$  &   2 & $x^6+a_2x^4+a_1x^2+1$ \\
 3&                        & $C_{7}$     &  7   & 1  & $7^4$  &   2 & $x^3+a_1x+a_2x^2+1$ \\
\hline
 4& \multirow{4}{*}{$D_{2m}$} & $D_4\times C_2$ &  2 & 2  & $2^8$     &  5 & $\prod_{i=1}^5(x^4+a_ix^2+1))$ \\
 4&                           & $ D_{10}\times C_2$  &  2 & 5  & $2^4$, 5     &  2 & $(x^{10}+a_1x^5+1)(x^{10}+a_2x^5+1)$ \\
 4&                           &$ D_{20}\times C_2$  &  2 & 10  & $2^3$, 10     &  1 & $x^{20}+a_1x^{10}+1$ \\
 4&                           &$ D_{4}\times C_4$  &  4 & 2  & $2^3, 4^2$      &  2 & $(x^4+a_1x^2+1)(x^4+a_2x^2+1)$ \\
 4&                           &$ D_{8}\times C_4$  &  4 & 4  & $2^2, 4^2$      &  1 & $x^8+a_1x^4+1$ \\
 5&                           &$G_5$  &  2 & 4  & $2^3$, $4^2$ &  2 & $(x^4-1)(x^8+a_1x^4+1)(x^8+a_2x^4+1)$  \\
 5&                           &$G_5$  &  2 & 20 & 2, 4, 20     &  0 & $x^{20}-1$  \\
 5&                           &$G_5$  &  4 & 8 & 2,$8^2$     &  0 & $x^{8}-1$  \\
 6&                           &$D_6\times C_2$  &  2 & 3  & $2^5$, 6    &  3 & $x \, \prod_{i=1}^3 (x^6+a_i x^3+1)$ \\
 6&                           &$D_{18}\times C_2$  &  2 & 9  & $2^3$, 18    &  1& $x(x^{18}+a_1x^9+1)$ \\
 6&                           &$D_{6}\times C_4$  &  4 & 3  & $2^2$, 4, 12    &  1& $x(x^6+a_1x^3+1)$ \\
 7&                           & $G_7$ &  2 & 2  & $2^5$, $4^2$ & 4 & $(x^4-1)\prod_{i=1}^4(x^4+a_ix^2+1)$ \\
 7&                           &$G_9$  &  2 & 5  & 2, $4^2$, 5 &  1 & $(x^{10}-1)(x^{10}+a_1x^5+1)$ \\
 7&                           & $G_7$ &  4 & 2  & 2, 4, $8^2$ & 1 & $(x^4-1)(x^4+a_1x^2+1)$ \\
 8&                           &$G_8$  &  2 & 2  & $2^5$, $4^2$ & 4 & $x(x^2-1)\prod_{i=1}^4(x^4+a_ix^2+1)$ \\
 8&                           & $G_8$ &  2 & 6 & $2^2$, 4, 12     &  1 & $x (x^{6} -1)(x^{12}+a_1x^6+1)$ \\
 8&                           & $G_8$ &  2 & 18 & 2, 4, 36     &  0 & $x (x^{18} -1)$ \\
 8&                           & $D_6\times C_3$ &  3 & 3 & 2, 3, 6, 9     &  1 & $x (x^{3} -1)(x^{6}+a_1x^3+1)$ \\
 8&                           & $D_{18}\times C_3$ &  3 & 9 & 2, 6, 27     & 0 & $x (x^{9}-1)$ \\
  8&                           & $G_8$ &  4 & 2 & 2, 4, $8^2$     &  1 & $x (x^{2} -1)(x^4+a_1x^2+1)$ \\
   8&                           & $G_8$ &  4 & 6 & 2, 8, 24     &  0 & $x (x^{6} -1)$ \\
  8&                           & $D_6\times C_7$ &  7 & 3 & 2, 14, 21     &  0 & $x (x^{3} -1)$ \\
   8&                           & $G_8$ &  10 & 2 & 2, $20^2$     &  0 & $x (x^{2} -1)$ \\
  9&                           & $G_9$ &  2 & 3 & $2^2, 4^2, 6$     &  2 & $x (x^{6} -1)(x^6+a_1x^3+1)(x^6+a_2x^3+1)$ \\
\hline
12 &   $A_4$                     & $K$ & 2 & 0 & $2^2$, $6^2$ & 1 & $(x^8+14x^4+1) \, f_1(x) $ \\
\hline
 17&   $S_4$                     &$G_{17}$  & 4  & 0 & 2, 4, 12 & 0 & $x^{8}+14x^4+1$ \\
 21&                             &$G_{21}$  & 2  & 0 & $4^2$ 6 & 0 & $(x^{8}+14x^4+1)(x^{12}-33x^8-33x^4+1)$ \\
\hline
27 &   $A_5$                     &  & 2 &  & 2, 5, 6 & 0 & $x^{20}-228x^{15}+494x^{10}+228x^5+1$ \\
\hline \hline
\multicolumn{8}{c}{Genus 10} \\
\hline \hline
 1& \multirow{4}{*}{$C_m$} &   $C_2^2$        &  2  & 2  & $2^{12}$      & 10 & $x^{22} + \sum_{i=1}^{10} a_i x^{2i} + 1$  \\
 1&                        &   $C_2 \times C_3$ &  3  & 2 & $2^2, 3^5$ & 5 & $x^{12} + \sum_{i=1}^5 a_i x^{2i} + 1$ \\
 1&                        &   $C_3^2$  &  3  & 3 & $3^5$ & 3  & $x^{12} + a_1 x^3 + a_2 x^{6}+a_3x^{9} +1$ \\
 1&                        &   $C_3\times C_4$  &  3  & 4 & $3^2, 4^2$ & 2  & $x^{12} + a_1 x^4 + a_2 x^{8} +1$ \\
 1&                        &   $C_2\times C_6$  &  6  & 2 & $2^2, 6^2$ & 2  & $x^{6} + a_1 x^2 + a_2 x^{4} +1$ \\
 2&                        & $C_{6}$  &  2  & 3 &   $2^6$, 3, 6   & 6   & $x^{21}+ \sum_{i=1}^6 a_i x^{3i} +1$ \\
 2&                        & $C_{14}$  &  2 & 7  &   $2^2$, 7, 14  &  2   & $x^{21} + a_1 x^7+a_2x^{14} +1$ \\
 2&                        & $C_{42}$  & 2 & 21  &  21, 42  &  0    & $x^{21} +1$ \\
 2&                        & $C_{33}$  &  3 & 11  &  11, 33   &  0    & $x^{11} +1$ \\
 2&                        & $C_{10}$  &  5 & 2  &   2, $5^2$, 10  &  2   & $x^6+a_2x^4+a_1x^2+1$ \\
 2&                        & $C_{15}$  &  5 & 3  &  3, 5, 15   & 1    & $x^6+a_1x^3 +1$ \\
 2&                        & $C_{30}$  &  5 & 6  &  6, 30    &  0      & $x^6 +1$ \\
\hline \hline
\end{tabular}
\end{table}

\addtocounter{table}{-1}
\begin{table}
\caption{(Cont.)}
\begin{tabular}{|l|l|l|l|l|l|l|l|}
\hline
Nr. & $\bar G$             & G&$n$  &$m$ & sig. & $\delta$ & Equation $y^n=f(x)$ \\
\hline 
 2&                        & $C_{30}$  &  6 & 5  &   5, 30   &  0      & $x^5  +1$ \\
 2&                        & $C_{33}$     &  11  & 3  &  3, 33  & 0 & $x^3 +1$ \\
 2&                        & $C_{42}$     &  21  & 2  & 2, 42       &   0 & $x^2+1$ \\
 3&                        & $C_2$     &  2  & 1  & $2^{21}$        &   19 & $x^{20}+\sum_{i=1}^{19} a_i x^{i} +1$ \\
 3&                        & $C_4$     &  2  & 2  & $2^9, 4^2$        &   9 & $x^{20}+\sum_{i=1}^9 a_i x^{2i} +1$ \\
 3&                        & $C_8$     &  2  & 4  & $2^4, 8^2$        &   4 & $x^{20}+a_1x^4+a_2x^8+a_3x^{12}+a_4x^{16} +1$ \\
 3&                        & $C_{10}$     &  2  & 5  & $2^3, {10}^2$  &   3 & $x^{20}+ a_1x^5+a_2x^{10}+a_3x^{15} +1$ \\
 3&                        & $C_{3}$     &  3   & 1  & $3^{11}$  &   9 & $x^{10}+ \sum_{i=1}^9 a_i x^{i}+1$ \\
 3&                        & $C_{6}$     &  3   & 2  & $3^4, 6^2$  &   4 & $x^{10}+ a_1x^2+a_2x^4+a_3x^6+a_4x^8+1$ \\
 3&                        & $C_{5}$     &  5   & 1  & $5^6$  &   4 & $x^{5}+ \sum_{i=1}^4 a_i x^{i}+1$ \\
 3&                        & $C_{6}$     & 6   & 1  & $6^5$  &   3 & $x^4+a_1x+a_2x^2a_3x^3+1$ \\
\hline
 4& \multirow{4}{*}{$D_{2m}$} & $D_{22}\times C_2$ &  2 & 11  & $2^3$, 11     &  1 & $x^{22}+a_1x^{11}+1$ \\
 4&                           & $ D_{4}\times C_3$  &  3 & 2  & $2^3, 3^3$     &  3 & $\prod_{i=1}^3 (x^4+a_i x^2+1) $ \\
 4&                           &$ D_{6}\times C_3$  &  3 & 3  & $2^2, 3^3$,    &  2 & $(x^6+a_1x^3+1)(x^6+a_2x^3+1)$ \\
 4&                           &$ D_{12}\times C_3$  &  3 & 6  & $2^2, 3, 6$      &  1 & $(x^{12}+a_1x^6+1$ \\
 4&                           &$ D_{6}\times C_6$  &  6 & 3  & $2^2, 3, 6$      &  1 & $x^6+a_1x^3+1$ \\
 5&                           &$G_5$  &  2 & 2  & $2^7$, 4 &  5 & $(x^2-1)\prod_{i=1}^5(x^4+a_ix^2+1)$  \\
 5&                           &$G_5$  &  2 & 22 & 2, 4, 22     &  0 & $x^{22}-1$  \\
 5&                           &$D_8\times C_3$  &  3 & 4 & 2, 3, 4, 6     &  1 & $(x^4-1)(x^8+a_1x^4+1)$  \\
 5&                           &$D_{24}\times C_3$  &  3 & 12  & 2, 6, 12   &  0& $x^{12}-1$ \\
 5&                           &$G_5$  &  6 & 2 & $2^2$, 6, 12     &  1 & $(x^{2}-1)(x^4+a_1x^2+1)$  \\
 5&                           &$G_5$  &  6 & 6 &    2, 6, 12     &  0 & $x^6-1$  \\
  6&                           &$D_8$  &  2 & 2 & $2^7$, 4     &  5 & $x\prod_{i=1}^5(x^4+a_ix^2+1)$  \\
 6&                           &$D_{10}\times C_2$  &  2 & 5  & $2^4$, 10    &  2& $x(x^{10}+a_1x^5+1)(x^{10}+a_2x^5+1)$ \\
 6&                           &$D_{40}$  &  2 & 10  & $2^3$, 20    &  1 & $x(x^{20}+a_1x^{10}+1)$ \\
 6&                           &$D_{10}\times C_3$  &  3 & 5  & $2^2$, 3, 15    &  1& $x(x^{10}+a_1x^5+1)$ \\
 6&                           &$D_{24}$  &  6 & 2  & $2^2$, 6, 12    &  1 & $x(x^{4}+a_1x^{2}+1)$ \\
 7&                           & $D_4\times C_3$ &  3 & 2  & 2, $3^2$, $6^2$ & 2 & $(x^2-1)(x^4+a_1x^2+1)(x^4+a_2x^2+1)$ \\
 7&                           &$D_6\times C_3$  &  3 & 3  &  $3^2, 6^2$ &  1 & $(x^{6}-1)(x^{6}+a_1x^3+1)$ \\
 8&                           &$G_8$  &  2 & 4  & $2^3$, 4, 8 & 2 & $x(x^4-1)(x^8+a_1x^4+1)(x^8+a_2x^4+1)$ \\
 8&                           & $G_8$ &  2 & 20 & 2, 4, 40    &  0 & $x (x^{20} -1)$ \\
 8&                           & $D_4\times C_3$ &  3 & 2 & 2,$3^2, 6^2$   &  2& $x(x^2-1)(x^4+a_1x^2+1)(x^4+a_2x^2+1)$ \\
 8&                           & $D_{20}\times C_3$ &  3 & 10 & 2, 6, 30   & 0 & $x (x^{10} -1)$ \\
 8&                           & $D_{10}\times C_5$ &  5 & 5 & 2, 10, 25   & 0 & $x (x^{5}-1)$ \\
 8&                           & $G_8$ &  6 & 4 & 2, 12, 24     &  0 & $x (x^{4} -1)$ \\
 8&                           & $D_4\times C_{11}$ &  11 & 2 & 2, $22^2$     &  0 & $x (x^{2} -1)$ \\
 9&                           & $G_9$ &  2 & 2 & $2^4, 4^3$     & 4 & $x(x^4-1)\prod_{i=1}^4(x^4+a_ix^2+1)$ \\
 9&                           & $G_9$ &  2 & 5 & $2, 4^2, 10$     & 1 & $x(x^{10}-1)(x^{10}+a_1x^5+1)$ \\
\hline
10 &   $A_4$                     & & 3 & 0 & 2, $3^3$ & 1 & $ f_1 (x) $ \\
14&                              & & 2 & 0 & 2, 3, 4, 6 & 1 & $x(x^4-1)(x^4+ 2 \sqrt{-3} \, x^2+1)\, f_1(x)$\\
\hline
 18&   $S_4$                     &$G_{18}$  & 6  & 0 & 2, 3, 24 & 0 & $x(x^4-1)$ \\
 20&                             &$S_4\times C_3$  & 3  & 0 & 3, 4, 6 & 0 & $x^{12}-33x^8-33x^4+1$ \\
\hline
25 &   $A_5$                     &$A_5\times C_3$  & 3 & 0 & 2, 3, 15 & 0 & $x(x^{10}+11x^5-1)$ \\
\hline
\end{tabular}
\end{table}

\clearpage

\section{Final remarks}

Following the methods described above we intend to create a database of all superelliptic curves of genus $g\leq 48$, inclusions among the loci, and the corresponding parametric equations for each family. Each locus can be determined in terms of the invariants of binary forms, but this is a difficult task computationally since such forms are not known for high degree binary forms.  However, such loci can also be described in terms of the dihedral invariants of superelliptic curves.

%*********************************************************************************************************
\nocite{*}
\bibliographystyle{plain}

\bibliography{myrefs}
%****************************************************************

\end{document}